\documentclass[11pt,a4paper]{article}
\usepackage[utf8]{inputenc}
\usepackage[T1]{fontenc}
\usepackage{amsmath, amssymb, amsthm}
\usepackage{geometry}
\usepackage{graphicx}
\usepackage{cite}
\usepackage{float}
\usepackage{placeins}

\begin{document}

\title{Oscillatory Interference in Dirichlet L-Functions and the Separation of Primes}

\author{Jouni J. Takalo}

\maketitle

\begin{abstract}

Dirichlet's theorem guarantees infinitely many primes in each reduced residue
class modulo $q$, but the analytic mechanism underlying this separation is
often difficult to visualize directly. In this article we construct simplified
oscillatory reconstructions based on the imaginary parts of the nontrivial
zeros of Dirichlet $L$-functions. 

These reconstructions produce interference patterns that act as analytic
filters separating primes according to congruence classes. Examples for
moduli $3$, $4$, and $5$ illustrate how the oscillatory frequencies associated
with the zeros generate structured peak patterns at prime powers. For complex
characters modulo $5$, conjugate pairs of $L$-functions produce cancellation
effects that mirror algebraic relations between characters. When all
characters modulo $5$ are combined, the Dedekind factorization of the
cyclotomic field $\mathbf{Q}(\zeta_5)$ appears visually as a striking
interference pattern in which only primes congruent to $1 \pmod 5$ remain.

These numerical experiments provide a visual bridge between analytic number
theory and algebraic number theory by illustrating how the zero distributions
of $L$-functions generate structured oscillations in prime-related functions.

\end{abstract}

\section{Introduction}

Dirichlet's theorem on arithmetic progressions states that each reduced
residue class modulo $q$ contains infinitely many primes. The analytic
framework behind this result is provided by Dirichlet characters and their
associated $L$-functions, which form one of the central tools of analytic
number theory \cite{Davenport_2000,IwaniecKowalski_2004}.

The explicit formula connects the distribution of prime numbers to the
nontrivial zeros of these $L$-functions. Through this formula, the zeros
generate oscillatory contributions to weighted prime counting functions.
Although the explicit formula is well known in analytic number theory,
the oscillatory structure produced by the zeros is often difficult to
visualize directly \cite{Edwards1974,MazurStein_2016}.

The purpose of this article is to illustrate this oscillatory structure
numerically. Using lists of imaginary parts of zeros of Dirichlet
$L$-functions, we construct simplified cosine--sine superpositions that
approximate the oscillatory terms appearing in the explicit formula.
While these reconstructions omit certain slowly varying weights, they
retain the essential interference patterns produced by the zeros.

These oscillatory reconstructions reveal a striking phenomenon: the
frequencies generated by the zeros act as analytic filters that separate
prime numbers according to congruence classes. Examples for moduli $3$,
$4$, and $5$ illustrate how different Dirichlet characters produce
distinct filtering patterns.

In the case of modulus $5$, the interaction between real and complex
characters reveals an especially clear interference structure. When the
oscillatory contributions corresponding to all characters modulo $5$ are
combined, the Dedekind factorization of the cyclotomic field
$\mathbf{Q}(\zeta_5)$ appears visually as a cancellation pattern in which
all residue classes except $p \equiv 1 \pmod 5$ disappear.

Modern computational tools make it possible to explore these phenomena
numerically. In this article the zero data were obtained using SageMath  \cite{Sagemath}
and the L-functions and Modular Forms Database (LMFDB)
\cite{LMFDB_2021}. The goal of the computations is not high-precision
evaluation but rather visualization of the oscillatory interference
patterns generated by the zeros.

\section*{Computational Setup}

The numerical reconstructions in this article are based on
lists of imaginary parts $\gamma$ of the nontrivial zeros
$\rho=\tfrac12+i\gamma$ of the relevant Dirichlet $L$-functions.

The zeros were computed using SageMath,
which provides built-in functionality for numerical evaluation
of Dirichlet $L$-functions and their zeros.
In each example, a large collection of consecutive zeros
was used in the simplified oscillatory reconstruction
described earlier.

The purpose of the computations here
is not high-precision evaluation,
but visual illustration of the oscillatory
interference patterns generated by the zeros.

\section*{The Explicit Formula (Conceptual View)}

For a Dirichlet character $\chi$ modulo $q$,
define the weighted prime counting function

\[
\Psi_\chi(x)=\sum_{n\le x}\chi(n)\Lambda(n),
\]

where $\Lambda(n)$ is the von Mangoldt function.

The explicit formula connects $\Psi_\chi(x)$
to the zeros of the associated $L$-function.
For a primitive character $\chi$,

\[
\Psi_\chi(x)
=
\text{Main Term}
-
\sum_{\rho} \frac{x^\rho}{\rho}
+
\text{(small correction terms)},
\]

where the sum runs over the nontrivial zeros
$\rho=\beta+i\gamma$ of $L(s,\chi)$.

If $\chi$ is trivial,
$L(s,\chi)$ has a pole at $s=1$,
and the main term is essentially $x$.
For nontrivial characters,
there is no pole at $s=1$,
and thus no dominant linear term.
In that case the oscillatory sum over zeros
governs the visible behavior.

The trivial zeros contribute only small smooth corrections.
The visible spikes in numerical reconstructions
arise from the nontrivial zeros.
For an accessible introduction to the explicit formula
and its relation to prime counting functions,
see Mazur and Stein~\cite{Edwards_2004, MazurStein_2016}. For a classical and more detailed treatment,
see Davenport~\cite{Davenport_2000}.

\section*{From Zeros to Oscillatory Waves}

If $\rho=\tfrac12+i\gamma$ is a nontrivial zero,
then

\[
x^\rho
=
x^{1/2+i\gamma}
=
x^{1/2}e^{i\gamma\log x}.
\]

Using Euler’s formula,

\[
e^{i\theta}=\cos\theta+i\sin\theta,
\]

each zero generates cosine and sine waves
in the variable $\log x$.

The full explicit formula contains additional weights
coming from the factor $1/\rho$
and from $x^{1/2}$.
However, to highlight the interference phenomenon clearly,
we use the following simplified oscillatory model.

\begin{center}
\fbox{
\parbox{0.9\linewidth}{
\textbf{Simplified Oscillatory Reconstruction}

\smallskip

\textbf{Input:}
Imaginary parts $\gamma$ of the zeros
$\rho=\tfrac12+i\gamma$.

\smallskip

\textbf{Construction:}
For each $x>0$,
compute
\[
S_\chi(x)
=
-\sum_{\gamma}
\cos(\gamma\log x),
\quad
T_\chi(x)
=
-\sum_{\gamma}
\sin(\gamma\log x).
\]

\smallskip

\textbf{Interpretation:}
Constructive interference of these waves
produces visible spikes at prime powers.
}
}
\end{center}

This simplified model omits the slowly varying weights
of the explicit formula,
but preserves the oscillatory structure
responsible for prime filtering.

\section*{Primes Modulo 3}

Modulo $3$ there are two characters:
the trivial character and one nontrivial real character.
The nontrivial character assigns

\[
\chi(n)=
\begin{cases}
1, & n\equiv 1 \pmod{3},\\
-1, & n\equiv 2 \pmod{3},\\
0, & 3\mid n.
\end{cases}
\]

Using the simplified oscillatory reconstruction,
we observe that primes congruent to $1$ modulo $3$
and primes congruent to $2$ modulo $3$
appear with opposite signs.
The zero frequencies act as an analytic filter
separating the two residue classes.

\begin{figure}[H]
\centering
\includegraphics[width=0.85\textwidth]{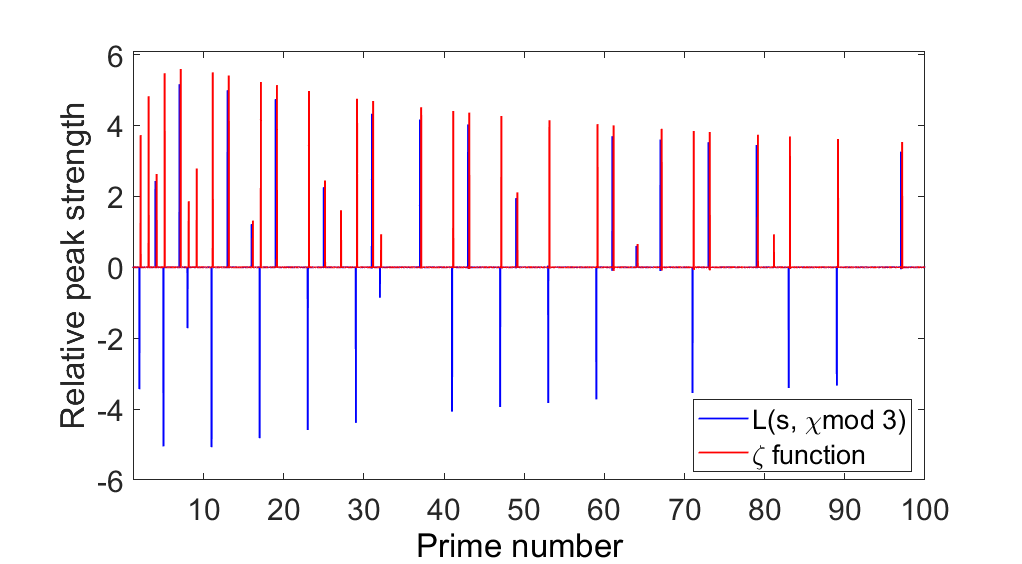}
\caption{
Oscillatory reconstruction for the nontrivial character
modulo $3$ (blue), compared with the zeta reconstruction (red).
The character separates the residue classes
$1$ and $2$ modulo $3$ into peaks of opposite sign.
}
\label{fig:mod3}
\end{figure}

\section*{Primes Modulo 4}

The multiplicative group $(\mathbb Z/4\mathbb Z)^\times=\{1,3\}$
admits one nontrivial character,

\[
\chi(n)=
\begin{cases}
1, & n\equiv 1 \pmod{4},\\
-1, & n\equiv 3 \pmod{4},\\
0, & 2\mid n.
\end{cases}
\]

The oscillatory reconstruction again separates
the two residue classes.
Primes congruent to $1$ modulo $4$
produce positive spikes,
while those congruent to $3$ modulo $4$
produce negative spikes.

This example is classical,
as primes congruent to $1$ modulo $4$
are precisely those expressible
as sums of two squares.

Since the character modulo $4$ is real,
the sine series vanishes identically,
and the reconstruction is purely real-valued.

\begin{figure}[H]
\centering
\includegraphics[width=0.85\textwidth]{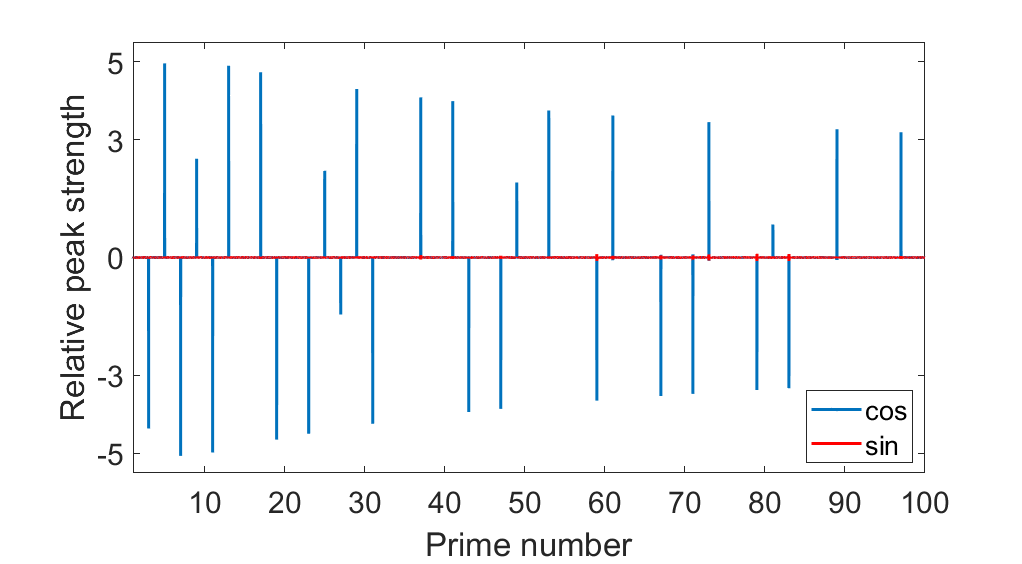}
\caption{
Oscillatory reconstruction for the nontrivial
Dirichlet character modulo $4$.
Since the character is real,
the sine contribution vanishes,
and the reconstruction is purely real.
Primes $p\equiv 1 \pmod{4}$ and $p\equiv 3 \pmod{4}$
appear as peaks of opposite sign.
}
\label{fig:mod4}
\end{figure}

\section*{Primes Modulo 5}

Modulo $5$ there are four Dirichlet characters:
one trivial character $\chi_0$,
one real quadratic character $\chi_2$,
and two complex conjugate characters $\chi_1$ and $\chi_3$.
Together they illustrate how real and complex characters
produce distinct filtering effects.

\subsection*{The Real Quadratic Character}

The quadratic character modulo $5$ takes values
\[
\chi_2(n)=
\begin{cases}
\;\;1, & n\equiv 1,4 \pmod{5},\\
-1, & n\equiv 2,3 \pmod{5},\\
\;\;0, & 5\mid n.
\end{cases}
\]

In the oscillatory reconstruction,
prime powers $p^k$ with $p\equiv 1,4 \pmod{5}$
appear as positive peaks,
while those with $p\equiv 2,3 \pmod{5}$
appear as negative peaks.
Prime powers of $5$ do not contribute.

Thus the quadratic character separates
quadratic residues from nonresidues modulo $5$.
The cosine series alone suffices,
since the character is real.

\begin{figure}[H]
\centering
\includegraphics[width=0.85\textwidth]{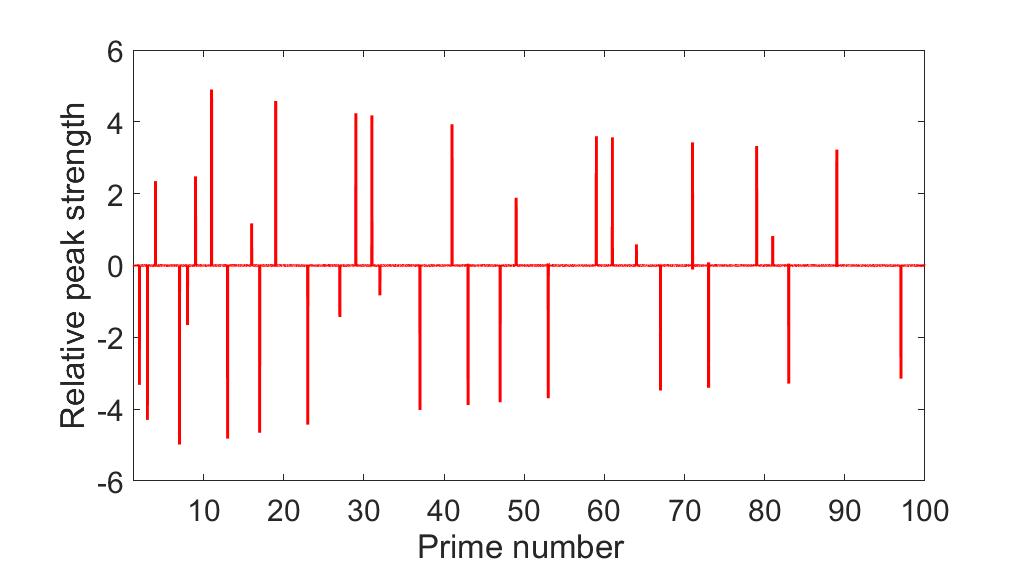}
\caption{
Oscillatory reconstruction for the quadratic
Dirichlet character modulo $5$.
Quadratic residues $1$ and $4$ modulo $5$
produce positive peaks,
while nonresidues $2$ and $3$
produce negative peaks.
}
\label{fig:mod5_quadratic}
\end{figure}

\subsection*{The Complex Conjugate Characters}

The characters $\chi_1$ and $\chi_3$
take values among the fourth roots of unity
$\{1,i,-1,-i\}$ and satisfy $\chi_3=\overline{\chi_1}$.

Their oscillatory reconstructions are complex-valued.
Examining the real and imaginary parts separately
reveals complementary behavior.

For $p\equiv 1,4 \pmod{5}$,
the character values are real,
and the peaks appear in the real part of the reconstruction.

For $p\equiv 2,3 \pmod{5}$,
the character values are purely imaginary,
and the peaks appear in the sine series.
Because $\chi_3$ is the complex conjugate of $\chi_1$,
their imaginary parts are equal in magnitude
but opposite in sign.
Consequently, when the two characters are combined,
their imaginary contributions cancel.

This cancellation illustrates how algebraic conjugation
manifests itself as destructive interference
in the oscillatory model.

\begin{figure}[H]
\centering
\includegraphics[width=0.9\textwidth]{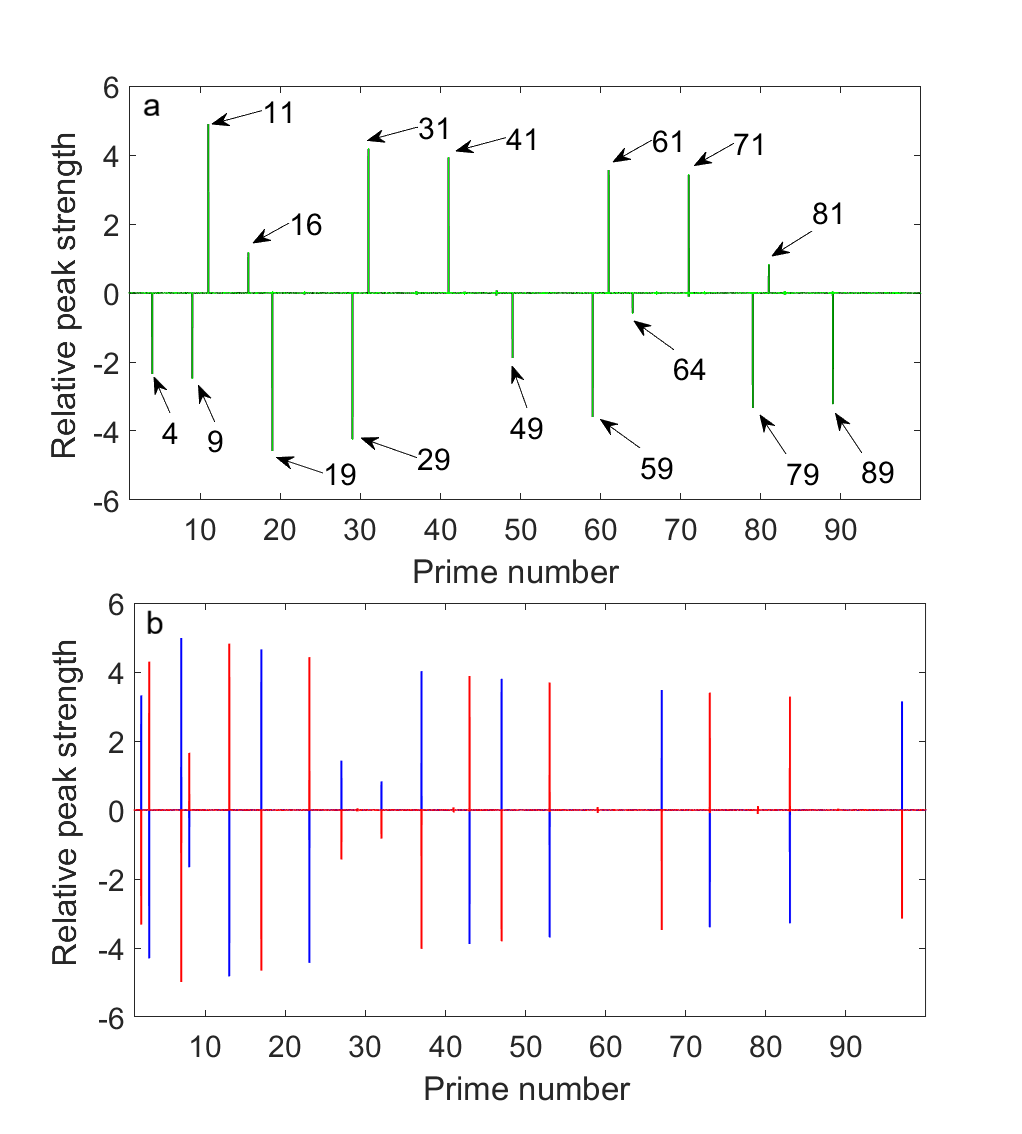}
\caption{
Oscillatory reconstructions for the complex conjugate
Dirichlet characters modulo $5$.
(a) Real parts of $L(s,\chi_1)$ and $L(s,\chi_3)$
produce identical peak patterns.
(b) Imaginary parts (blue for $\chi_1$, red for $\chi_3$)
have equal magnitude but opposite sign,
illustrating cancellation between conjugate characters.
}
\label{fig:mod5_complex}
\end{figure}

\subsection*{Dedekind Zeta Factorization}

The deeper algebraic structure emerges
when all characters modulo $5$ are combined.

The Dedekind zeta function of the cyclotomic field
$\mathbb{Q}(\zeta_5)$ factorizes as
\[
\zeta_{\mathbb{Q}(\zeta_5)}(s)
=
\prod_{\chi \bmod 5} L(s,\chi)
=
\zeta(s)\,
L(s,\chi_2)\,
L(s,\chi_1)\,
L(s,\chi_3).
\]

Its nontrivial zeros are precisely the union
of the zeros of these four factors.

When we sum the oscillatory contributions
coming from all four zero sets,
a striking simplification occurs.

Prime powers $p\equiv 2,3 \pmod{5}$
vanish by destructive interference:
their contributions are opposite in
$\zeta(s)$ and $L(s,\chi_2)$,
and zero in the real parts of
$L(s,\chi_1)$ and $L(s,\chi_3)$.

Prime powers $p\equiv 4 \pmod{5}$
also cancel,
since their contributions are positive
for $\zeta(s)$ and $L(s,\chi_2)$
but negative in the real parts of
$L(s,\chi_1)$ and $L(s,\chi_3)$.

Only the class $p\equiv 1 \pmod{5}$
survives constructively in all four factors.

Moreover, the imaginary parts from
$L(s,\chi_1)$ and $L(s,\chi_3)$
cancel identically,
so the resulting reconstruction is purely real.

Thus the Dedekind factorization,
an algebraic identity,
appears visually as a perfect interference pattern:
all residue classes cancel except $1 \pmod{5}$.

One exception must be noted.
The class $p\equiv 1 \pmod 5$ consists of primes
that split completely in $\mathbb{Q}(\zeta_5)$.
The remaining visible spikes at $5^k$
arise solely from the Riemann zeta factor,
reflecting the ramification of $5$.

\begin{figure}[H]
\centering
\includegraphics[width=0.9\textwidth]{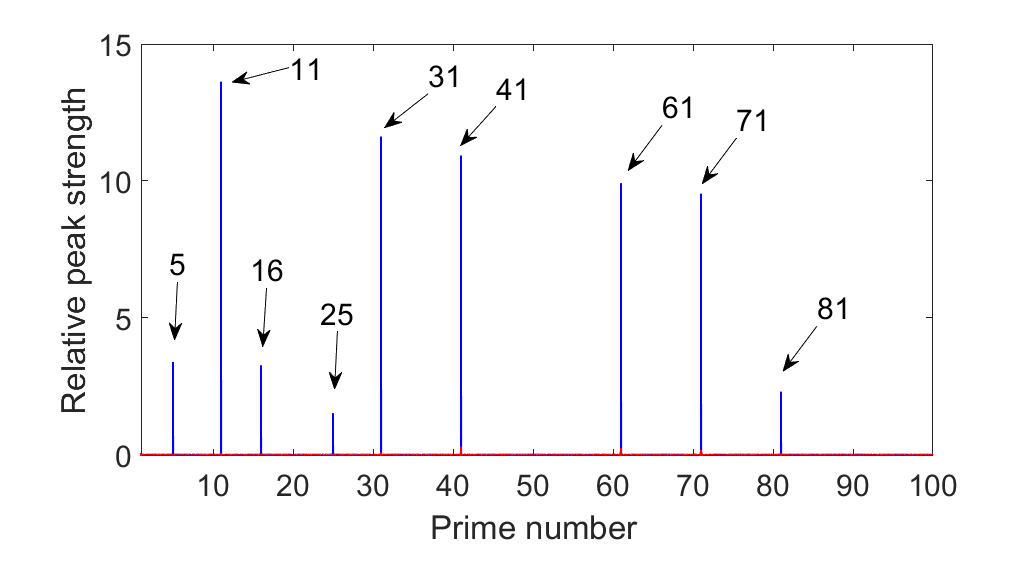}
\caption{
Reconstruction corresponding to the Dedekind zeta function
$\zeta_{\mathbb{Q}(\zeta_5)}(s)$,
obtained by combining all characters modulo $5$.
Destructive interference eliminates every residue class
except $p\equiv 1 \pmod{5}$,
while the prime powers $5^k$ persist due to ramification.
The imaginary part (red) vanishes.
}
\label{fig:dedekind}
\end{figure}

\section*{Why Modulo 6 Adds Nothing New}

The multiplicative group
\[
(\mathbb Z/6\mathbb Z)^\times = \{1,5\}
\]
is isomorphic to $(\mathbb Z/3\mathbb Z)^\times$.
Since modulo $2$ admits only the trivial character,
no new nontrivial Dirichlet characters arise modulo $6$.
Thus working modulo $6$ introduces no new oscillatory behavior
beyond what already appears in the case $q=3$.

\section*{Suggested Explorations}

Readers may wish to:

\begin{itemize}
\item Vary the number of zeros used in the reconstruction
and observe how the sharpness of spikes changes.
\item Plot the oscillatory sums as functions of $\log x$
to see more regular interference patterns.
\item Experiment with modulus $7$
and compare real and complex characters.
\end{itemize}

\section*{Conclusion}

The zeros of Dirichlet $L$-functions encode oscillatory
frequencies that act as analytic filters for prime numbers.
Even in the simplified cosine–sine model,
their interference separates primes into congruence classes.

For real characters, such as those modulo $3$ and $4$,
the reconstruction is purely real and visibly distinguishes
residue classes by sign.
For complex characters modulo $5$,
real and imaginary parts reveal complementary structures,
and conjugate pairs produce exact cancellation of
imaginary contributions.

The Dedekind factorization for $\mathbb{Q}(\zeta_5)$
provides a striking culmination:
when all four character contributions are combined,
destructive interference eliminates every residue class
except $p\equiv 1 \pmod{5}$.
An algebraic product identity thus appears visually
as a perfect interference pattern.

These numerical reconstructions offer a concrete bridge
between analytic and algebraic number theory,
showing how zero distributions manifest themselves
as structured oscillations in the primes.

\section*{Acknowledgments}
Zero data for L-functions was obtained using SageMath and Modular Forms Database (LMFDB).


\end{document}